\documentclass[11pt]{amsart}
\usepackage[all]{xy}
\usepackage{amsmath}
\usepackage{amsfonts}
\usepackage{amssymb}
\usepackage{latexsym}
\usepackage{graphicx}
\usepackage{color}
\usepackage{amscd}
\usepackage{epsfig}
\usepackage{cite}

\newtheorem{problem}{Problem}[section]
\newtheorem{definition}[problem]{Definition}
\newtheorem{lemma}[problem]{Lemma}
\newtheorem{theorem}[problem]{Theorem}

\newtheorem{corollary}[problem]{Corollary}

\newtheorem{conjecture}[problem]{Conjecture}
\newtheorem{example}{Example}[section]

\title{T-adic exponential sums over finite fields}
\author{Chunlei Liu}
\address{Department of Mathematical Sciences, Shanghai Jiao Tong
University, Shanghai 200240, P.R. China, E-mail: clliu@sjtu.edu.cn}
\author{Daqing Wan}
\address{Department of Mathematics, University of California,
Irvine, CA 92697-3875, USA, E-mail: dwan@math.uci.edu}
\begin{document}
\maketitle
\begin{abstract}
$T$-adic exponential sums associated to a Laurent polynomial $f$ are
introduced. They interpolate all classical $p^m$-power order
exponential sums associated to $f$. The Hodge bound for the Newton
polygon of $L$-functions of $T$-adic exponential sums is
established. This bound enables us to determine, for all $m$, the
Newton polygons of $L$-functions of $p^m$-power order exponential
sums associated to an $f$ which is ordinary for $m=1$. Deeper
properties of $L$-functions of $T$-adic exponential sums are also
studied. Along the way, new open problems about the $T$-adic
exponential sum itself are discussed.
\end{abstract}



\section{Introduction}
\subsection{Classical exponential sums}We first recall the definition of classical exponential sums over finite
fields of characteristic $p$ with values in a $p$-adic field.

Let $p$ be a fixed prime number, $\mathbb{Z}_p$ the ring of $p$-adic
integers, $\mathbb{Q}_p$ the field of $p$-adic numbers, and
$\overline{\mathbb{Q}}_p$ a fixed algebraic closure of
$\mathbb{Q}_p$. Let $q=p^a$ be a power of $p$, $\mathbb{F}_q$ the
finite field of $q$ elements, $\mathbb{Q}_q$ the unramified
extension of $\mathbb{Q}_p$ with residue field $\mathbb{F}_q$, and
$\mathbb{Z}_q$ the ring of integers of $\mathbb{Q}_q$.

Fix a positive integer $n$. Let
$f(x)\in\mathbb{Z}_q[x_1^{\pm1},x_2^{\pm1},\cdots,x_n^{\pm1}]$ be a
Laurent polynomial in $n$ variables  of the form
$$f(x)=\sum\limits_{u}a_ux^u, \ a_u\in\mu_{q-1},\ x^u =x_1^{u_1}\cdots x_n^{u_n},$$ where $\mu_k$
denotes the group of $k$-th roots of unity in
$\overline{\mathbb{Q}}_p$.

\begin{definition}Let $\psi$ be a locally constant character of
$\mathbb{Z}_p$ of order $p^m$ with values in
$\overline{\mathbb{Q}}_p$, and let $\pi_{\psi}=\psi(1)-1$. The sum
$$S_{f,\psi}(k)= \sum\limits_{x\in\mu_{q^k-1}^n}
\psi(\text{Tr}_{\mathbb{Q}_{q^k}/\mathbb{Q}_p}(f(x)))$$ is called a
$p^m$-power order exponential sum on the $n$-torus ${\Bbb G}_m^n$
over $\mathbb{F}_{q^k}$. The generating function
$$L_{f,\psi}(s) =L_{f,\psi}(s;\mathbb{F}_q)
=\exp(\sum\limits_{k=1}^{\infty}S_{f,\psi}(k)\frac{s^k}{k})\in
1+s\mathbb{Z}_p[\pi_{\psi}][[s]]$$ is called the $L$-function of
$p^m$-power order exponential sums over $\mathbb{F}_q$ associated to
$f(x)$.
\end{definition}

Note that the above exponential sum for $m\geq 1$ is still an
exponential sum over a finite field as we just sum over the subset
of roots of unity (corresponding to the elements of a finite field
via the Teichm\"uller lifting), not over the whole finite residue
ring $\mathbb{Z}_q/p^m\mathbb{Z}_q$. The exponential sum over the
whole finite ring $\mathbb{Z}_q/p^m\mathbb{Z}_q$ and its generating
function as $m$ varies is the subject of Igusa's zeta function, see
Igusa \cite{Ig}.

In general, the above $L$-function $L_{f,\psi}(s)$ of exponential
sums is rational in $s$. But, if $f$ is non-degenerate, then
$L_{f,\psi}(s)^{(-1)^{n-1}}$ is a polynomial, as was shown in
\cite{AS, AS2} for $\psi$ of order $p$, and in \cite{LW} for all
$\psi$. By a result of \cite{GKZ}, if $p$ is large enough, then $f$
is generically non-degenerate. For non-degenerate $f$, the location
of the zeros of $L_{f,\psi}(s)^{(-1)^{n-1}}$ becomes an important
issue. The $p$-adic theory of such L-functions was developed by
Dwork, Bombieri \cite{B}, Adolphson-Sperber \cite{AS, AS2}, the
second author \cite{Wa1, Wa2}, and Blache \cite{B4} for $\psi$ of
order $p$. More recently initial part of the theory was extended to
all $\psi$ by Liu-Wei \cite{LW} and Liu \cite{Liu}.

The $p$-adic theory of the above exponential sum for $n=1$ and
$\psi$ of order $p$  has a long history and has been studied
extensively in the literature. For instance, in the simplest case
 that $f(x)=x^d$, the exponential
sum was studied by Gauss, see Berndt-Evans \cite{BE} for a
comprehensive survey. By the Hasse-Davenport relation for Gauss
sums, the L-function is a polynomial whose zeros are given by roots
of Gauss sums. Thus, the slopes of the L-function are completely
determined by the Stickelberger theorem for Gauss sums. The roots of
the L-function have explicit $p$-adic formulas in terms of $p$-adic
$\Gamma$-function via the Gross-Koblitz formula \cite{GK}. These
ideas can be extended to treat the so-called diagonal $f$ case for
general $n$, see Wan \cite{Wa2}. These elementary cases have been
used as building bricks to study the deeper non-diagonal $f(x)$ via
various decomposition theorems, which are the main ideas of Wan
\cite{Wa1,Wa2}. In the case $n=1$ and $\psi$ of order $p$, further
progresses about the slopes of the L-function were made in Zhu
\cite{Zh1,Zh2}, Blache and Ferard \cite{BF}, and Liu \cite{Liu3}.

\subsection{$T$-adic exponential sums}We now define the $T$-adic exponential
sum, state our main results, and put forward some new questions.

\begin{definition}For a positive integer $k$, the $T$-adic exponential sum of $f$ over $\mathbb{F}_{q^k}$
is the sum:
$$S_f(k,T)=\sum\limits_{x\in\mu_{q^k-1}^n}
(1+T)^{\text{\rm
Tr}_{\mathbb{Q}_{q^k}/\mathbb{Q}_p}(f(x))}\in\mathbb{Z}_p[[T]].$$
The $T$-adic $L$-function of $f$ over $\mathbb{F}_q$ is the
generating function
$$L_f(s,T) =L_f(s,T;\mathbb{F}_q) =\exp(\sum\limits_{k=1}^{\infty}S_f(k,T)\frac{s^k}{k})\in 1+s\mathbb{Z}_p[[T]][[s]].$$
\end{definition}

The $T$-adic exponential sum interpolates classical exponential sums
of $p^m$-order over finite fields for all positive integers $m$. In
fact, we have
$$S_f(k,\pi_{\psi})=S_{f,\psi}(k).$$ Similarly, one can recover
the classical L-function of the $p^m$-order exponential sum from the
$T$-adic $L$-function by the formula
$$L_f(s,\pi_{\psi}) =L_{f,\psi}(s).$$

We view $L_f(s, T)$ as a power series in the single variable $s$
with coefficients in the complete discrete valuation ring
$\mathbb{Q}_p[[T]]$ with uniformizer $T$.

\begin{definition}The $T$-adic characteristic function of $f$ over $\mathbb{F}_q$, or $C$-function of $f$ for short, is
the generating function $$C_f(s,T)
=\exp(\sum\limits_{k=1}^{\infty}-(q^k-1)^{-n}S_{f}(k,T)\frac{s^k}{k})\in
1+s\mathbb{Z}_p[[T]][[s]].$$
\end{definition}

The $C$-function $C_f(s,T)$ and the $L$-function $L_f(s,T)$
determine each other. They are related by
$$L_f(s,T) = \prod_{i=0}^n C_f(q^is, T)^{(-1)^{n-i-1}{n\choose
i}},$$and
$$C_f(s, T)^{(-1)^{n-1}}= \prod_{j=0}^{\infty} L_f(q^js, T)^{n+j-1\choose j}.$$

In \S4, we prove
\begin{theorem}[analytic continuation] The
$C$-function $C_f(s,T)$ is $T$-adic entire in $s$. As a consequence,
the $L$-function $L_f(s,T)$ is $T$-adic meromorphic in $s$.
\end{theorem}
The above theorem tells that the $C$-function behaves $T$-adically
better than the $L$-function. In fact, in the $T$-adic setting, the
$C$-function is a more natural object than the L-function. Thus, we
shall focus more on the $C$-function.

Knowing the analytic continuation of $C_f(s,T)$, we are then
interested in the location of its zeros. More precisely, we would
like to determine the $T$-adic Newton polygon of this entire
function $C_f(s,T)$. This is expected to be a complicated problem in
general. It is open even in the simplest case $n=1$ and $f(x)=x^d$
is a monomial if $p\not\equiv 1~(\mod d)$. What we can do is to give
an explicit combinatorial lower bound  depending only on $q$ and
$\Delta$, called the $q$-Hodge bound ${\rm HP}_q(\Delta)$. This
polygon will be described in detail in \S3.

Let ${\rm {NP}}_T(f)$ denote the $T$-adic Newton polygon of the
$C$-function $C_f(s,T)$.
 In
\S5, we prove
\begin{theorem}[Hodge bound] We have
$$ {\rm {NP}}_T(f) ~\geq ~{\rm HP}_q(\Delta).$$
\end{theorem}

This theorem shall give several new results on classical exponential
sums, as we shall see in \S2. In particular, this extends, in one
stroke,  all known ordinariness results for $\psi$ of order $p$ to
all $\psi$ of any $p$-power order. It demonstrates the significance
of the $T$-adic L-function. It also gives rise to the following
definition.
\begin{definition}The Laurent polynomial $f$ is called $T$-adically ordinary
if ${\rm NP}_{T}(f)~=~{\rm HP}_q(\Delta)$.
\end{definition}

We shall show that the classical notion of ordinariness implies
$T$-adic ordinariness. But it is possible that a non-ordinary $f$ is
$T$-adically ordinary. Thus, it remains of interest to study exactly
when $f$ is $T$-adically ordinary. For this purpose, in \S6, we
extend the facial decomposition theorem in Wan \cite{Wa1} to the
$T$-adic case. Let $\Delta$ be the convex closure in $\mathbb{R}^n$
of the origin and the exponents of the non-zero monomials in the
Laurent polynomial $f(x)$. For any closed face $\sigma$ of $\Delta$,
we let $f_{\sigma}$ denote the sum of monomials of $f$ whose
exponent vectors lie in $\sigma$.

\begin{theorem}[$T$-adic facial decomposition]\label{facial-decomposition}
The Laurent polynomial $f$ is $T$-adically ordinary if and only if for every closed face
$\sigma$ of $\Delta$ of codimension $1$ not containing the origin,
the restriction $f_{\sigma}$ is $T$-adically ordinary.
\end{theorem}

In \S7, we briefly discuss the variation of the $C$-function
$C_f(s,T)$ and its Newton polygon when the reduction of $f$ moves in
an algebraic family over a finite field. The main questions are the
generic ordinariness, generic Newton polygon, the analogue of the
Adolphson-Sperber conjecture \cite{AS}, Wan's limiting conjecture
\cite{Wa2}, Dwork's unit root conjecture \cite{Dw0} in the $T$-adic
and $\pi_{\psi}$-adic case. We shall give an overview about what can
be proved and what is unknown, including a number of conjectures.
Basically, a lot can be proved in the ordinary case, and a lot
remain to be proved in the non-ordinary case.

{\bf Acknowledgement. }The first author is supported by NSFC Grant
No. 10671015.

\section{Applications}
In this section we give several applications of the $T$-adic
exponential sum to classical exponential sums.

\begin{theorem}[integrality theorem]We have$$L_f(s,T) \in 1+s\mathbb{Z}_p[[T]][[s]],$$
and
$$C_f(s,T) \in 1+s\mathbb{Z}_p[[T]][[s]].$$
\end{theorem}
\proof Let $|\mathbb{G}_m^n|$ be the set of closed points of
$\mathbb{G}_m^n$ over $\mathbb{F}_q$, and $a\mapsto\hat{a}$ the
Teichm\"{u}ller lifting. It is easy to check that the $T$-adic
L-function has the Euler product expansion
$$L_f(s,T)=\prod\limits_{x\in |\mathbb{G}_m^n|}{1\over
(1-(1+T)^{\text{Tr}_{\mathbb{Q}_{q^{\deg(x)}}/\mathbb{Q}_p}(f(\hat{x}))}s^{\deg(x)})}\in
1+s\mathbb{Z}_p[[T]][[s]],$$ where
$\hat{x}=(\hat{x}_1,\cdots,\hat{x}_n)$. The theorem now
follows.\qed\hfill

The above proof shows that the L-function $L_f(s,T)$ is the
L-function $L(s,\rho_f)$ of the following continuous $(p,T)$-adic
representation of the arithmetic fundamental group:
$$\rho_f: \pi_1^{\rm arith}({\mathbb G}_m^n/\mathbb{F}_q)
\longrightarrow {\rm GL}_1(\mathbb{Z}_p[[T]]),$$ defined by
$$\rho_f({\rm Frob}_x)
=(1+T)^{\text{Tr}_{\mathbb{Q}_{q^{\deg(x)}}/\mathbb{Q}_p}(f(\hat{x}))}.$$
The rank one representation $\rho_f$ is transcendental in nature.
Its L-function $L(s,\rho_f)$ seems to be beyond the reach of
$\ell$-adic cohomology, where $\ell$ is a prime different from $p$.
However, the specialization of $\rho_f$ at the special point
$T=\pi_{\psi}$ is a character of finite order. Thus, the
specialization
$$L(s,\rho_f)|_{T=\pi_{\psi}} =L_{f, \psi}(s)$$
can indeed be studied using Grothendieck's $\ell$-adic trace formula
\cite{Gr}. This gives another proof that the L-function $L_{f,
\psi}(s)$ is a rational function in $s$. But the $T$-adic L-function
$L_f(s,T)$ itself is certainly out of the reach of $\ell$-adic
cohomology as it is truly transcendental.

Let ${\rm {NP}}_T(f)$ denote the $T$-adic Newton polygon of the
$C$-function $C_f(s,T)$, and let ${\rm {NP}}_{\pi_{\psi}}(f)$ denote
the $\pi_{\psi}$-adic Newton polygon of the $C$-function
$C_f(s,\pi_{\psi})$. The integrality of $C_f(s,T)$ immediately gives
the following theorem.
\begin{theorem}[rigidity bound] If $\psi$ is non-trivial, then$$ {\rm {NP}}_{\pi_{\psi}}(f) ~\geq ~{\rm {NP}}_T(f).$$
\end{theorem}
\proof Obvious.\qed\hfill

A natural question is to ask when ${\rm {NP}}_{\pi_{\psi}}(f)$
coincides with its rigidity bound.
\begin{theorem}[transfer theorem]If $ {\rm {NP}}_{\pi_{\psi}}(f) ~= ~{\rm {NP}}_T(f)$ holds for one non-trivial
$\psi$, then it holds for all non-trivial $\psi$.
\end{theorem}
\proof By the integrality of $C_f(s,T)$, the $T$-adic Newton polygon
of $C_f(s,T)$ coincides with the $\pi_{\psi}$-adic Newton polygon of
$C_f(s,\pi_{\psi})$ if and only if for every vertex $(i,e)$ of the
$T$-adic Newton polygon of $C_f(s,T)$, the coefficients of $s^i$ in
$C_f(s,T)$ differs from $T^e$ by a unit in
$\mathbb{Z}_p[[T]]^{\times}$. It follows that if the coincidence
happens for one non-trivial $\psi$, it happens for all non-trivial
$\psi$. The theorem is proved.\qed\hfill

\begin{definition}We call $f$ rigid if $ {\rm {NP}}_{\pi_{\psi}}(f) ~= ~{\rm {NP}}_T(f)$ for one (and hence for all)
non-trivial $\psi$.
\end{definition}

In \cite{LLN}, cooperating with his students, the first author
showed that $f$ is generically rigid if $n=1$ and $p$ is
sufficiently large. So the rigid bound is the best possible bound.
In contrast, the weaker Hodge bound ${\rm HP}_q(\Delta)$ is only
best possible if $p\equiv 1~(\mod ~d)$, where $d$ is the degree of
$f$.

We now pause to describe the relationship between the Newton
polygons of $C_f(s,\pi_{\psi})$ and $L_{f,\psi}(s)^{(-1)^{n-1}}$. We
need the following definitions.
\begin{definition}A convex polygon with initial point $(0,0)$ is
called algebraic if it is the graph of a $\mathbb{Q}$-valued
function defined on $\mathbb{N}$ or on an interval of $\mathbb{N}$,
and its slopes are of finite multiplicity and of bounded
denominator.
\end{definition}
\begin{definition}For an algebraic polygon with slopes $\{\lambda_i\}$, we define
its slope series to be $\sum\limits_it^{\lambda_i}$.
\end{definition}

It is clear that an algebraic polygon is uniquely determined by its
slope series. So the slope series embeds the set of algebraic
polygons into the ring
$\lim\limits_{\stackrel{\rightarrow}{d}}\mathbb{Z}[[t^{\frac{1}{d}}]]$.
The image is
$\lim\limits_{\stackrel{\rightarrow}{d}}\mathbb{N}[[t^{\frac{1}{d}}]]$.
It is closed under addition and multiplication. Therefore one can
define an addition and a multiplication on the set of algebraic
polygons.

\begin{lemma}\label{a-l-relation}Suppose that $f$ is non-degenerate.
Then the $q$-adic Newton polygon of $C_f(s,\pi_{\psi};\mathbb{F}_q)$
is the product of the $q$-adic Newton polygon of
$L_{f,\psi}(s;\mathbb{F}_q)^{(-1)^{n-1}}$ and the algebraic polygon
$\frac{1}{(1-t)^n}$.
\end{lemma}

\proof Note that the $C$-value $C_f(s,\pi_{\psi})$ and the
$L$-function $L_{f,\psi}(s)$ determine each other. They are related
by
$$L_{f,\psi}(s) = \prod_{i=0}^n C_f(q^is, \pi_{\psi})^{(-1)^{n-i-1}{n\choose
i}},$$and
$$C_f(s, \pi_{\psi})^{(-1)^{n-1}}= \prod_{j=0}^{\infty} L_{f,\psi}(q^js)^{n+j-1\choose j}.$$
Suppose that
$$L_{f,\psi}(s)^{(-1)^{n-1}}=\prod_{i=1}^d(1-\alpha_is).$$
Then
$$C_f(s,\pi_{\psi}) =
\prod_{j=0}^{\infty}\prod_{i=1}^d(1-\alpha_iq^{j}s)^{n+j-1\choose
j}.$$ Let $\lambda_i$ be the $q$-adic order of $\alpha_i$. Then the
$q$-adic order of $\alpha_iq^{j}$ is $\lambda_i+j$. So the slope
series of the $q$-adic Newton polygon of
$L_{f,\psi}(s)^{(-1)^{n-1}}$ is
$$S(t)=\sum\limits_{i=1}^dt^{\lambda_i},$$ and the slope series
of the $q$-adic Newton polygon of $C_f(s,\pi_{\psi})$ is
$$\sum\limits_{j=0}^{+\infty}\sum\limits_{i=0}^d {n+j-1\choose
j}t^{\lambda_i+j}=\frac{1}{(1-t)^n}S(t).$$ The lemma now
follows.\qed\hfill

We combine the rigidity bound and the Hodge bound to give the
following theorem.
\begin{theorem}If $\psi$ is non-trivial, then
$$ {\rm {NP}}_{\pi_{\psi}}(f) ~\geq ~{\rm {NP}}_T(f) ~\geq ~{\rm HP}_q(\Delta).$$
\end{theorem}
\proof Obvious.\qed\hfill

If we drop the middle term, we arrive at the Hodge bound $${\rm
{NP}}_{\pi_{\psi}}(f)~\geq ~{\rm HP}_q(\Delta)$$ of
Adolphson-Sperber \cite{AS2} and Liu-Wei \cite{LW}.

\begin{theorem}If $ {\rm {NP}}_{\pi_{\psi}}(f) ~= ~{\rm HP}_q(\Delta)$
holds for one non-trivial $\psi$, then $f$ is rigid, $T$-adically
ordinary, and the equality holds for all non-trivial $\psi$.
\end{theorem}
\proof Suppose that ${\rm {NP}}_{\pi_{\psi_0}}(f) ~= ~{\rm
HP}_q(\Delta)$ for a non-trivial $\psi_0$. Then, by the last
theorem, we have$${\rm {NP}}_{\pi_{\psi_0}}(f) ~= ~{\rm {NP}}_T(f)
~= ~{\rm HP}_q(\Delta).$$ So $f$ is rigid and $T$-adically ordinary,
and $${\rm {NP}}_{\pi_{\psi}}(f) ~= ~{\rm {NP}}_T(f) ~= ~{\rm
HP}_q(\Delta)$$ holds for all nontrivial $\psi$. The theorem is
proved.\qed\hfill

\begin{definition}We call $f$ ordinary
if ${\rm NP}_{\pi_{\psi}}(f)~=~{\rm HP}_q(\Delta)$ holds for one
(and hence for all) non-trivial $\psi$.
\end{definition}

The notion of ordinariness now carries much more information than
what we had known. From this, we see that the $T$-adic exponential
sum provides a new framework to study all $p^m$-power order
exponential sums simultaneously. Instead of the usual way of
extending the methods for $\psi$ of order $p$ to the case of higher
order, the $T$-adic exponential sum has the novel feature that it
can sometimes transfer a known result for one non-trivial $\psi$ to
all non-trivial $\psi$. This philosophy is carried out further in
the paper \cite{LLN}.

\begin{example}Let $$f(x)=x_1+x_2+\cdots+x_n+\frac{\alpha}{x_1x_2\cdots x_n}, \
\alpha\in\mu_{q-1}.$$ Then, by the result of Sperber \cite{Sp} and
our new information on ordinariness, we have
$${\rm
NP}_{\pi_{\psi}}(f)~=~{\rm HP}_q(\Delta)$$ for all non-trivial
$\psi$.\end{example}
\section{The $q$-Hodge polygon}

In this section, we describe explicitly the $q$-Hodge polygon mentioned in the introduction.
Recall that $f(x)\in\mathbb{Z}_q[x_1^{\pm1},x_2^{\pm1},\cdots,x_n^{\pm1}]$
is a Laurent polynomial in $n$ variables  of the form
$$f(x)=\sum\limits_{u\in
\mathbb{Z}^n}a_ux^u, \ a_u\in {\mathbb{Z}_q}, \ a_u^q=a_u.$$
We stress that the non-zero coefficients of $f(x)$ are roots of unity
in $\mathbb{Z}_q$, thus correspond in a unique way to
Teichm\"uller liftings of elements of the finite field $\mathbb{F}_q$.
If the coefficients of $f(x)$ are arbitrary elements in $\mathbb{Z}_q$,
much of the theory still holds, but it is more complicated to describe
the results. We have made the simplifying assumption that the
non-zero coefficients are always roots of unity in this paper.

Let $\Delta$ be
the convex polyhedron in $\mathbb{R}^n$ associated to $f$, which is
generated by the origin and the exponent vectors of the non-zero
 monomials of $f$. Let $C(\Delta)$ be the cone in $\mathbb{R}^n$
 generated by $\Delta$. Define the degree function $u\mapsto\deg(u)$
 on $C(\Delta)$ such that $\deg(u)=1$ when $u$ lies on a
 codimensional $1$ face of $\Delta$ that does not contain the origin,
 and such that
 $$\deg(ru)=r\deg(u), \ r\in \mathbb{R}_{\geq 0}, \ u \in
 C(\Delta).$$ We call it the degree function associated to $\Delta$.
 We have $\deg(u+v)\leq\deg(u)+\deg(v)$ if $u,v\in C(\Delta)$, and
 the equality holds if and only if $u$ and $v$ are co-facial. In
 other words, the number
 $$c(u,v):=\deg(u)+\deg(v)-\deg(u+v)$$
 is $0$ if $u,v\in C(\Delta)$ are co-facial, and is positive
 otherwise. We call that number $c(u,v)$ the co-facial defect of $u$
 and $v$. Let
 $$M(\Delta):=C(\Delta)\cap\mathbb{Z}^n$$ be the
 set of lattice points in the cone $C(\Delta)$. Let $D$ be the
 denominator of the degree function, which is the smallest positive
 integer such that
 $$\deg M(\Delta)\subset\frac{1}{D}\mathbb{Z}.$$
 For every natural number $k$, we define
 $$W(k):=W_{\Delta}(k)= \#\{ u\in M(\Delta) | \deg(u)=k/D\}$$ to be
 the number of lattice points of degree $\frac{k}{D}$ in $M(\Delta)$.
For prime power $q=p^a$,
the $q$-Hodge polygon of $f$ is the polygon with vertices $(0,0)$
 and
 $$(\sum\limits_{j=0}^iW(j),
 a(p-1)\sum\limits_{j=0}^i\frac{j}{D}W(j)),\ i=0,1,\cdots.
 $$
 It is also called the $q$-Hodge polygon of $\Delta$ and denoted by
 ${\rm HP}_q(\Delta)$. It depends only on $q$ and $\Delta$.
 It has a side of slope $a(p-1){\frac{j}{D}}$ with horizontal length
 $W(j)$ for each non-negative integer $j$.

\section{Analytic continuation}

 In this section, we prove the $T$-adic analytic continuation of the
 C-function $C_f(s,T)$. The idea is to employ Dwork's trace formula
 in the $T$-adic case.

 Note that the Galois group $\text{Gal}(\mathbb{Q}_q/\mathbb{Q}_p)$
 is cyclic of order $a=\log_pq$. There is an element in the Galois
 group whose restriction to $\mu_{q-1}$ is the $p$-power morphism. It
 is of order $a$, and is called the Frobenius element. We denote that
 element by $\sigma$.

 We define a new variable $\pi$ by the relation $E(\pi)=1+T$, where
 $$E(\pi)=\exp(\sum_{i=0}^{\infty}\frac{\pi^{p^i}}{p^i}) \in 1+\pi{\mathbb Z}_p[[\pi]]$$
 is the Artin-Hasse exponential series. Thus, $\pi$ and $T$ are two
 different uniformizers of the $T$-adic local ring ${\mathbb Q}_p[[T]]$.
 It is clear that for $\alpha\in \mathbb{Z}_q$, we have
 $$E(\pi\alpha)\in 1+\pi{\mathbb Z}_q[[\pi]],$$
 and for $\beta\in \mathbb{Z}_p$, we have
 $$E(\pi)^{\beta}\in 1+\pi{\mathbb Z}_p[[\pi]].$$

 The Galois group $\text{Gal}(\mathbb{Q}_q/\mathbb{Q}_p)$ can act on
 $\mathbb{Z}_q[[\pi]]$ but keeping $\pi$ fixed. The Artin-Hasse
 exponential series has a kind of commutativity expressed as the
 following lemma.

 \begin{lemma}[Commutativity]\label{commutativity}

 We have the following commutative diagram
$$
 \begin{array}{ccc}
 \mu_{q-1} & \stackrel{E(\pi\cdot)}{\rightarrow} & \mathbb{Z}_q[[\pi]] \\
 \text{\rm Tr}\downarrow &  & \downarrow\text{\rm Norm}\\
 \mu_{p-1} & \stackrel{E(\pi)^{\cdot}}{\rightarrow} & \mathbb{Z}_p[[\pi]]. \\
 \end{array}
$$
That is, if $x\in\mu_{q-1}$, then
$$E(\pi)^{x+x^p+\cdots+x^{p^{a-1}}}=E(\pi x)E(\pi x^p)\cdots E(\pi x^{p^{a-1}}).$$
\end{lemma} {\it Proof.  }Since for $x\in\mu_{q-1}$,
$$\sum\limits_{j=0}^{a-1}x^{p^j}=\sum\limits_{j=0}^{a-1}x^{p^{j+i}},$$
we have
$$
E(\pi)^{x+x^p+\cdots+x^{p^{a-1}}}
=\exp(\sum_{i=0}^{\infty}\frac{\pi^{p^i}}{p^i}\sum\limits_{j=0}^{a-1}x^{p^{j+i}})=E(\pi
x)E(\pi x^p)\cdots E(\pi x^{p^{a-1}}).$$ The lemma is proved.
\hfill\qed

\begin{definition}
Let $\pi^{1/D}$ be a fixed $D$-th root of $\pi$. Define
$$L(\Delta)=\{\sum_{u\in M(\Delta)}b_u\pi^{\deg(u)}x^u:\
 b_u\in\mathbb{Z}_q[[\pi^{1/D}]] \},$$
 and
$$B=\{\sum\limits_{u\in M(\Delta)}b_u\pi^{\deg(u)}x^u \in L(\Delta),\
 \text{\rm ord}_T(b_u)\rightarrow+\infty
 \text{ if }\deg(u)\rightarrow+\infty\}.$$
\end{definition}

The spaces $L(\Delta)$ and $B$ are $T$-adic Banach algebras over the
ring $\mathbb{Z}_{q}[[\pi^{\frac{1}{D}}]]$. The monomials
$\pi^{\deg(u)}x^u$ ($u\in M(\Delta))$ form an orthonormal basis
(resp., a formal basis) of $B$ (resp., $L(\Delta))$ . The algebra
$B$ is contained in the larger Banach algebra $L(\Delta)$. If $u\in
\Delta$, it is clear that $E(\pi x^u) \in L(\Delta)$. Write
$$E_f(x) :=\prod\limits_{a_u\neq0}E(\pi a_ux^u),\text{ if }f(x)=\sum\limits_{u\in \mathbb{Z}^n}a_ux^u.$$
This is an element of $L(\Delta)$ since $L(\Delta)$ is a ring.

The Galois group $\text{Gal}(\mathbb{Q}_q/\mathbb{Q}_p)$ can act on
$L(\Delta)$ but keeping $\pi^{1/D}$ as well as the variables $x_i$'s
fixed. From the commutativity of the Artin-Hasse exponential series,
one can infer the following lemma.
\begin{lemma}[Dwork's splitting lemma]If
 $x\in\mu_{q^k-1}$, then $$E(\pi)^{\text{\rm Tr}_{\mathbb{Q}_{q^k}/\mathbb{Q}_p}(f(x))}
 =\prod\limits_{i=0}^{ak-1}E_f^{\sigma^i}(x^{p^i})
,$$ where  $a$ is the order of $\text{\rm
Gal}(\mathbb{Q}_q/\mathbb{Q}_p)$.
\end{lemma} {\it Proof.  }We have
$$E(\pi)^{\text{Tr}_{\mathbb{Q}_{q^k}/\mathbb{Q}_p}(f(x))}
=\prod\limits_{a_u\neq0}E(\pi)^{\text{Tr}_{\mathbb{Q}_{q^k}/\mathbb{Q}_p}(a_ux^u)}$$
$$=\prod\limits_{a_u\neq0}\prod\limits_{i=0}^{ak-1}E(\pi(a_ux^u)^{p^i})=
\prod\limits_{i=0}^{ak-1}E_f^{\sigma^i}(x^{p^i}).$$The lemma is
proved. \hfill\qed

\begin{definition} We define a map
$$\psi_p:L(\Delta)\rightarrow L(\Delta),\
 \sum\limits_{u\in
M(\Delta)} b_ux^u\mapsto\sum\limits_{u\in M(\Delta)} b_{pu}x^u.$$
\end{definition}

It is clear that the composition map $\psi_p\circ E_f$ sends $B$ to
$B$.

\begin{lemma} Write
$$E_f(x) =\sum\limits_{u\in M(\Delta)}\alpha_u(f)\pi^{\deg(u)}x^u.$$ Then,
$\psi_p\circ E_f(\pi^{\deg(u)}x^u)$
$$=\sum\limits_{w\in
M(\Delta)}\alpha_{pw-u}(f)\pi^{c(pw-u,u)}\pi^{(p-1)\deg(w)}\pi^{\deg(w)}x^w,\
u\in M(\Delta),$$ where $c(pw-u, u)$ is the co-facial defect
of $pw-u$ and $u$.
\end{lemma} {\it Proof. } This follows directly from the definition
of $\psi_p$ and $E_f(x)$. \hfill\qed

\begin{definition} Define
$$\psi:=\sigma^{-1}\circ\psi_p\circ E_f: B\longrightarrow B,$$
and its $a$-th iterate
$$\psi^{a}=\psi_p^{a}\circ
\prod\limits_{i=0}^{a-1}E_{f}^{\sigma^i}(x^{p^i}).$$
\end{definition}

Note that $\psi$ is linear over $\mathbb{Z}_p[[\pi^{\frac{1}{D}}]]$,
but semi-linear over $\mathbb{Z}_q[[\pi^{\frac{1}{D}}]]$. On the
other hand, $\psi^a$ is linear over $\mathbb{Z}_q[[\pi^{1/D}]]$. By
the last lemma, $\psi^a$ is completely continuous in the sense of
Serre \cite{Se}.

\begin{theorem}[Dwork's trace formula]\label{trace-formula}For
every positive integer $k$,
$$(q^k-1)^{-n}S_f(k,T)
=\text{\rm
Tr}_{B/\mathbb{Z}_{q}[[\pi^{\frac{1}{D}}]]}(\psi^{ak}).$$\end{theorem}
{\it Proof. }Let $g(x)\in B$. We have
$$\psi^{ak}(g)=\psi_p^{ak}(g\prod\limits_{i=0}^{ak-1}E_f^{\sigma^i}(x^{p^i})).$$Write
$$\prod\limits_{i=0}^{ak-1}E_f^{\sigma^i}(x^{p^i})=\sum\limits_{u\in
M(\Delta)}\beta_u x^u.$$ One computes that
$$\psi^{ak}(\pi^{\deg(v)}x^v)=\sum\limits_{u\in
M(\Delta)}\beta_{q^ku-v}\pi^{\deg(v)}x^u.$$ Thus,
$${\rm Tr}(\psi^{ak}|B/\mathbb{Z}_{q}[[\pi^{\frac{1}{D}}]]) =
\sum\limits_{u\in M(\Delta)} \beta_{(q^k-1)u}.$$But, by Dwork's
splitting lemma, we have
$$(q^k-1)^{-n}S_f(k,T)=(q^k-1)^{-n}\sum\limits_{x\in\mu_{q^k-1}^n}
\prod\limits_{i=0}^{ak-1} E_f^{\sigma^i}(x^{p^i})=\sum\limits_{u\in
M(\Delta)} \beta_{(q^k-1)u}.$$ The theorem now follows. \hfill\qed

\begin{theorem}[Analytic trace formula]\label{analytic-trace-formula}
We have
$$C_f(s,T)=\det(1-\psi^as\mid
B/\mathbb{Z}_{q}[[\pi^{\frac{1}{D}}]]).$$ In particular, the
$T$-adic C-function $C_f(s,T)$ is $T$-adic analytic in $s$.
\end{theorem}

{\it Proof. }It follows from the last theorem and the well known
identity
$$\det(1-\psi^a s)
=\exp(-\sum\limits_{k=1}^{\infty}\text{Tr}(\psi^{ak})\frac{s^k}{k}).$$
\hfill\qed

This theorem gives another proof that the coefficients of $C_f(s,T)$
and $L_f(s,T)$ as power series in $s$ are $T$-adically integral.

\begin{corollary}For each non-trivial $\psi$, the $C$-value $C_f(s,\pi_{\psi})$ is
$p$-adic entire in $s$ and the  $L$-function $L_{f,\psi}(s)$ is
rational in $s$.
\end{corollary}
\proof Obvious.\qed\hfill
\section{The Hodge bound}

The analytic trace formula in the previous section reduces the study of $C_f(s,T)$ to the
study of the operator $\psi^a$. We consider $\psi$ first. Note that $\psi$
operates on $B$ and is linear over
$\mathbb{Z}_p[[\pi^{\frac{1}{D}}]]$.

\begin{theorem}\label{hodge-for-psi}The $T$-adic Newton polygon of $\det(1-\psi s\mid
B/\mathbb{Z}_p[[\pi^{\frac{1}{D}}]])$ lies above the polygon with
vertices $(0,0)$ and
$$(a\sum\limits_{k=0}^iW(k),
a(p-1)\sum\limits_{k=0}^i\frac{k}{D}W(k)),\
i=0,1,\cdots.$$\end{theorem}

{\it Proof. }Let $\xi_1,\xi_2,\cdots,\xi_a$ be a normal basis of
$\mathbb{Q}_q$ over $\mathbb{Q}_p$. Write
$$(\xi_j\alpha_{pw-u}(f))^{\sigma^{-1}}=\sum\limits_{i=0}^{a-1}\alpha_{(i,w),(j,u)}(f)\xi_i,\
\alpha_{(i,w),(j,u)}(f)\in\mathbb{Z}_p[[\pi^{1/D}]].$$ Then
$\psi(\xi_j\pi^{\deg(u)}x^u)$
$$=\sum\limits_{i=0}^{a-1}\sum\limits_{w\in
M(\Delta)}
\alpha_{(i,w),(j,u)}(f)\pi^{c(pw-u,u)}\pi^{(p-1)\deg(w)}\xi_i\pi^{\deg(w)}x^w.$$
That is, the matrix of $\psi$ over
$\mathbb{Z}_p[[\pi^{\frac{1}{D}}]]$ with respect to the orthonormal
basis $\{\xi_j\pi^{\deg(u)}x^u\}_{0\leq j<a,u\in M(\Delta)}$ is
$$A=(\alpha_{(i,w),(j,u)}(f)\pi^{c(pw-u,u)}\pi^{(p-1)\deg(w)})_{(i,w),(j,u)}.$$
So, the $T$-adic Newton polygon of $\det(1-\psi s\mid
B/\mathbb{Z}_p[[\pi^{\frac{1}{D}}]])$ lies above the polygon with
vertices $(0,0)$ and
$$(a\sum\limits_{k=0}^iW(k),
a(p-1)\sum\limits_{k=0}^i\frac{k}{D}W(k))\ (i=0,1,\cdots).$$ Theorem
\ref{hodge-for-psi} is proved. \hfill\qed

We are now ready to prove the Hodge bound for the Newton polygon.

\begin{theorem}\label{hodge-bound} We have
$${\rm {NP}}_T(f) ~\geq ~{\rm HP}_q(\Delta).$$
\end{theorem}

\proof By the above theorem, it suffices to prove that the $T$-adic
Newton polygon of
 $\det(1-\psi^as^a\mid
B/\mathbb{Z}_{q}[[\pi^{\frac{1}{D}}]])$ coincides with that of
$\det(1-\psi s\mid B/\mathbb{Z}_p[[\pi^{\frac{1}{D}}]])$. Note that
$$\det(1-\psi^a s\mid B/\mathbb{Z}_p[[\pi^{\frac{1}{D}}]])
={\rm Norm}(\det(1-\psi^a s\mid
B/\mathbb{Z}_q[[\pi^{\frac{1}{D}}]])),$$ where the norm map is the
norm from $\mathbb{Z}_q[[\pi^{\frac{1}{D}}]]$ to
$\mathbb{Z}_q[[\pi^{\frac{1}{D}}]]$. The theorem now follows from
the equality
$$\prod\limits_{\zeta^a=1}\det(1-\psi\zeta s\mid B/\mathbb{Z}_p[[\pi^{\frac{1}{D}}]])
=\det(1-\psi^a s^a\mid
B/\mathbb{Z}_p[[\pi^{\frac{1}{D}}]]).$$\hfill\qed

\section{Facial decomposition}

In this section, we extend the facial decomposition theorem in
\cite{Wa1}. Recall that the operator $\psi=\sigma^{-1}\circ
(\psi_p\circ E_f)$ is only semi-linear over
$\mathbb{Z}_q[[\pi^{\frac{1}{D}}]]$. But its second factor
$\psi_p\circ E_f$ is clearly linear and so det$(1-(\psi_p\circ E_f)
s\mid B/\mathbb{Z}_q[[\pi^{\frac{1}{D}}]])$ is well defined. We
begin with the following theorem.

\begin{theorem}\label{removal-of-sigma}The $T$-adic Newton polygon of
$f$ coincides with ${\rm HP}_q(\Delta)$ if and only if the $T$-adic
Newton polygon of det$(1-(\psi_p\circ E_f) s\mid
B/\mathbb{Z}_q[[\pi^{\frac{1}{D}}]]$) coincides with the polygon
with vertices $(0,0)$ and
$$(\sum\limits_{k=0}^iW(k),
(p-1)\sum\limits_{k=0}^i\frac{k}{D}W(k)),\
i=0,1,\cdots.$$\end{theorem}

{\it Proof. }In the proof of Theorem \ref{hodge-bound}, we showed
that the $T$-adic Newton polygon of
 $C_f(s^a,T)$
coincides with that of $\det(1-\psi s\mid
B/\mathbb{Z}_p[[\pi^{\frac{1}{D}}]])$. Note that
$$\det(1-(\psi_p\circ E_f) s\mid
B/\mathbb{Z}_p[[\pi^{\frac{1}{D}}]])={\rm Norm}(\det(1-(\psi_p\circ E_f) s\mid
B/\mathbb{Z}_{q}[[\pi^{\frac{1}{D}}]])),$$
where the norm map is the norm from $\mathbb{Z}_q[[\pi^{\frac{1}{D}}]]$
to $\mathbb{Z}_q[[\pi^{\frac{1}{D}}]]$.
The theorem is equivalent
to the statement that the $T$-adic Newton polygon of $\det(1-\psi
s\mid B/\mathbb{Z}_p[[\pi^{\frac{1}{D}}]])$ coincides with the
polygon with vertices $(0,0)$ and
$$(\sum\limits_{k=0}^iaW(k),
a(p-1)\sum\limits_{k=0}^i\frac{k}{D}W(k)),\ i=0,1,\cdots$$ if and
only if the $T$-adic Newton polygon of $\det(1-(\psi_p\circ E_f)
s\mid B/\mathbb{Z}_p[[\pi^{\frac{1}{D}}]])$  does. Therefore it
suffices to show that the determinant of the matrix
$$(\alpha_{(i,w),(j,u)}(f)\pi^{c(pw-u,u)})_{0\leq
i,j<a,\deg(w),\deg(u)\leq\frac{k}{D}}$$is not divisible by $T$ in
$\mathbb{Z}_q[[\pi^{\frac{1}{D}}]]$ if and only if the determinant
of the matrix
$$(\alpha_{pw-u}(f)\pi^{c(pw-u,u)})_{\deg(w),\deg(u)\leq\frac{k}{D}}$$
 is not divisible by $T$ in
$\mathbb{Z}_q[[\pi^{\frac{1}{D}}]]$. The theorem now follows from
the fact that the latter determinant is the norm of the former from
$\mathbb{Q}_q[[\pi^{\frac{1}{D}}]]$ to
$\mathbb{Q}_p[[\pi^{\frac{1}{D}}]]$ up to a sign. \hfill\qed

We now define the open facial decomposition $F(\Delta)$. It is the
decomposition of $C(\Delta)$ into a disjoint union of relatively
open cones generated by the relatively open faces of $\Delta$ whose
closure does not contain the origin. Note that every relatively open
cone generated by co-facial vectors in $C(\Delta)$ is contained in a
unique element of $F(\Delta)$.
\begin{lemma}\label{restriction-estimate}Let $\sigma\in F(\Delta)$, and $u\in\sigma$.
Then $\alpha_u(f_{\bar{\sigma}})\equiv\alpha_u(f)(\mod\pi^{1/D})$,
where $f_{\bar{\sigma}}$ is the sum of monomials of $f$ whose
exponent vectors lie in the closure $\bar{\sigma}$ of $\sigma$.
\end{lemma}
{\it Proof. }Let $v_1,\cdots,v_j$ be exponent vectors of monomials
of $f$ such that $a_1v_1+\cdots+a_jv_j=u$  with $a_1>0,\cdots,
a_j>0$. It suffices to show that either $v_1,\cdots v_j$ lie in the
closure of $\sigma$, or their contribution to $\alpha_u(f)$ is
$\equiv0(\mod\pi^{1/D})$. Suppose that their contribution to
$\alpha_u(f)$ is $\not\equiv0(\mod\pi^{1/D})$. Then $v_1,\cdots,v_j$
must be co-facial. So the interior of the cone generated by those
vectors is contained in a unique element of $F(\Delta)$. As that
interior has a common point $u$ with $\sigma$, it must be $\sigma$.
It follows that $v_1,\cdots v_j$ lie in the closure of $\sigma$. The
lemma is proved. \hfill\qed

\begin{lemma}Let $\sigma,\tau\in F(\Delta)$ be distinct. Let $w\in\sigma$, and $u\in\tau$.
Suppose that the dimension of $\sigma$ is no greater than that of
$\tau$. Then $pw-u$ and $u$ are not co-facial, i.e., $c(pw-u,u)>0$.
\end{lemma}
{\it Proof. }Suppose that $pw-u$ and $u$ are co-facial. Then the
interior of the cone generated by $pw-u$ and $u$ is contained in a
unique element of $F(\Delta)$. As that interior has a common point
$w$ with $\sigma$, it must be $\sigma$. It follows that $u$ lies in
the closure of $\sigma$. As $\sigma$ and $\tau$ are distinct, $u$
lies in the boundary of $\sigma$. This implies that the dimension of
$\tau$ is less than that of $\sigma$, which is a contradiction.
Therefore $pw-u$ and $u$ are not co-facial. The lemma is proved.
\hfill\qed

For $\sigma \in
F(\Delta)$, we define
$$M({\sigma}) =M(\Delta)\cap \sigma =\mathbb{Z}^n\cap \sigma$$
be the set of lattice points in the cone $\sigma$.

\begin{theorem}[Open facial decomposition]\label{open-facial-decomposition}The
$T$-adic Newton polygon of $f$ coincides with ${\rm HP}_q(\Delta)$
if and only if for every $\sigma\in F(\Delta)$, the determinants of
the matrices
$$\{\alpha_{pw-u}(f_{\bar\sigma})\pi^{c(pw-u,u)}\}_{w,u\in
M({\sigma}),\deg(w),\deg(u)\leq \frac{k}{D}},\ k=0,1,\cdots$$ are
not divisible by $T$ in $\mathbb{Z}_q[[\pi^{\frac{1}{D}}]]$, where
$\bar{\sigma}$ is the closure of $\sigma$.\end{theorem}

{\it Proof. }By Theorem \ref{removal-of-sigma}, the $T$-adic Newton
polygon of $C_f(s,T)$ coincides with the $q$-Hodge polygon of $f$ if and
only if the determinants of the matrices
$$A^{(k)}=\{\alpha_{pw-u}(f)\pi^{c(pw-u,u)}\}_{w,u\in
M(\Delta),\deg(w),\deg(u)\leq \frac{k}{D}},\ k=0,1,\cdots$$ are not
divisible by $T$ in $\mathbb{Z}_q[[\pi^{\frac{1}{D}}]]$. Write
$$A_{\sigma,\tau}^{(k)}=\{\alpha_{pw-u}(f)\pi^{c(pw-u,u)}\}_{w\in
M({\sigma}),u\in M({\tau}),\deg(w),\deg(u)\leq \frac{k}{D}}.$$ The
facial decomposition shows that $A^{(k)}$ has the block form
$(A_{\sigma,\tau}^{(k)})_{\sigma,\tau\in F(\Delta)}$. The last lemma
shows that the block form modulo $\pi^{\frac{1}{D}}$ is triangular
if we order the cones in $F(\Delta)$ in dimension-increasing order.
It follows that $\det A^{(k)}$ is not divisible by $T$ in
$\mathbb{Z}_q[[\pi^{\frac{1}{D}}]]$ if and only if for all
$\sigma\in F(\Delta)$, $\det A_{\sigma,\sigma}^{(k)}$ is not
divisible by $T$ in $\mathbb{Z}_q[[\pi^{\frac{1}{D}}]]$. By Lemma
\ref{restriction-estimate}, modulo $\pi^{\frac{1}{D}}$,
$A_{\sigma,\sigma}^{(k)}$ is congruent to the matrix
$$\{\alpha_{pw-u}(f_{\bar\sigma})\pi^{c(pw-u,u)}\}_{w,u\in
M({\sigma}),\deg(w),\deg(u)\leq \frac{k}{D}}.$$ So $\det
A_{\sigma,\sigma}^{(k)}$ is not divisible by $T$ in
$\mathbb{Z}_q[[\pi^{\frac{1}{D}}]]$ if and only if the determinant
of the matrix
$$\{\alpha_{pw-u}(f_{\bar\sigma})\pi^{c(pw-u,u)}\}_{w,u\in
M({\sigma}),\deg(w),\deg(u)\leq \frac{k}{D}}$$ is not divisible by
$T$ in $\mathbb{Z}_q[[\pi^{\frac{1}{D}}]]$.  The theorem is proved.
\hfill\qed

The closed facial decomposition Theorem \ref{facial-decomposition}
follows from the open decomposition theorem and the fact that
$$F(\Delta)=\bigcup_{\sigma\in
F(\Delta),\dim\sigma=\dim\Delta}F(\bar{\sigma}).$$ A similar
$\pi_{\psi}$-adic facial decomposition theorem for $C_f(s,
\pi_{\psi})$ can be proved in a similar way. Alternatively, it
follows from the transfer theorem together with the
$\pi_{\psi}$-adic facial decomposition in \cite{Wa1} for $\psi$ of
order $p$.

\section{Variation of C-functions in a family}

Fix an $n$-dimensional integral convex polytope $\triangle$ in
$\mathbb{R}^n$ containing the origin. For each prime $p$, let
$P(\Delta, \mathbb{F}_p)$ denote the parameter space of all Laurent
polynomials $f(x)$ over $\bar{\mathbb{F}}_p$ such that
$\Delta(f)=\Delta$. This is a connected rational variety defined
over $\mathbb{F}_p$. For each $f\in P(\Delta,
\mathbb{F}_p)(\mathbb{F}_q)$, the Teichm\"uller lifting
gives a Laurent polynomial $\tilde{f}$ whose non-zero coefficients are
roots of unity in $\mathbb{Z}_q$. The
C-function $C_{\tilde{f}}(s,T)$ is then defined and $T$-adically entire.
For simplicity of notation, we shall just write $C_f(s,T)$ for $C_{\tilde{f}}(s,T)$,
similarly, $L_f(s,T)$ for $L_{\tilde{f}}(s,T)$. Thus, our C-function and L-function are now
defined for Laurent polynomials over finite fields, via the Teichm\"uller lifting.
We would like to study how $C_f(s,T)$ varies when $f$ varies in the
algebraic variety $P(\Delta, \mathbb{F}_p)$.

Recall that for a closed face $\sigma\in \Delta$, $f_{\sigma}$
denotes the restriction of $f$ to $\sigma$. That is, $f_{\sigma}$ is
the sum of those non-zero monomials in $f$ whose exponents are in
$\sigma$.

\begin{definition} A Laurent polynomial $f\in P(\Delta, \mathbb{F}_p)$ is called
non-degenerate if for every closed face $\sigma$ of $\Delta$ of
arbitrary dimension which does not contain the origin, the system
$$\frac{\partial f_{\sigma}}{\partial x_1} =\cdots =\frac{\partial f_{\sigma}}{\partial
x_n} =0$$ has no common zeros with $x_1\cdots x_n\not=0$ over the
algebraic closure of $\mathbb{F}_p$.
\end{definition}

The non-degenerate condition is a geometric condition which insures
that the associated Dwork cohomology can be calculated. In
particular, it implies that, if $\psi$ is of order $p^m$, then the
L-function $L_{f,\psi}(s)^{(-1)^{n-1}}$ is a polynomial in $s$ whose
degree is precisely $n!{\rm Vol}(\Delta)p^{n(m-1)}$, see \cite{LW}.
As a consequence, we deduce
\begin{theorem}Let $f\in P(\Delta, \mathbb{F}_p)(\mathbb{F}_q)$. Write
$$L_f(s,T)^{(-1)^{n-1}} = \sum_{k=0}^{\infty} L_{f,k}(T)s^k,
\ L_{f,k}(T)\in \mathbb{Z}_p[[T]].$$ Assume that $f$ is
non-degenerate. Then for every positive integer $m$ and all positive
integer $k>n!{\rm Vol}(\Delta)p^{n(m-1)}$ , we have the following
congruence in $\mathbb{Z}_p[[T]]$:
$$L_{f,k}(T) \equiv 0~ \left(\mod \frac{(1+T)^{p^m}-1}{T}\right).$$
\end{theorem}

\proof  Write $$\frac{(1+T)^{p^m}-1}{T} =\prod (T-\xi).$$ The
non-degenerate assumption implies that
$$L_f(s,\xi)^{(-1)^{n-1}}=\sum_{j=0}^{\infty} L_{f,j}(\xi)s^j,$$ is a polynomials in $s$ of
degree $\leq n!{\rm Vol}(\Delta)p^{n(m-1)}<k$. It follows that
$L_{f,k}(\xi)=0$ for all $\xi$. That is, $L_{f,k}(T)$ is divisible
by $(T-\xi)$ for $\xi$. The theorem now follows. \hfill\qed

\begin{definition}Let $N(\Delta, \mathbb{F}_p)$ denote the subset of all
non-degenerate Laurent polynomials $f\in P(\Delta, \mathbb{F}_p)$.
\end{definition}

The subset $N(\Delta, \mathbb{F}_p)$ is Zariski open in $P(\Delta,
\mathbb{F}_p)$. It can be empty for some pair
$(\Delta,\mathbb{F}_p)$. But, for a given $\Delta$, $N(\Delta,
\mathbb{F}_p)$ is Zariski open dense in $P(\Delta, \mathbb{F}_p)$
for all primes $p$ except for possibly finitely many primes
depending on $\Delta$. It is an interesting and independent question
to classify the primes $p$ for which $N(\Delta, \mathbb{F}_p)$ is
non-empty. This is related to the GKZ discriminant \cite{GKZ}. For
simplicity, we shall only consider non-degenerate $f$ in the
following.

\subsection{Generic ordinariness} The first question is how often $f$ is $T$-adically
ordinary when $f$ varies in the non-degenerate locus $N(\Delta,
\mathbb{F}_p)$. Let $U_p(\Delta,T)$ be the subset of $f\in N(\Delta,
\mathbb{F}_p)$ such that $f$ is $T$-adically ordinary, and
$U_p(\Delta)$ the subset of $f\in N(\Delta, \mathbb{F}_p)$ such that
$f$ is ordinary. One can prove

\begin{lemma}The set
$U_p(\Delta)$ is Zariski open in $N(\Delta, \mathbb{F}_p)$.
\end{lemma}
One can ask if $U_p(\Delta, T)$ is also Zariski open in $N(\Delta,
\mathbb{F}_p)$. We do not know the answer.

Our question is for which $p$, $U_p(\Delta)$ and $U_p(\Delta,T)$ are
Zariski dense in $N(\Delta, \mathbb{F}_p)$. The rigidity bound as
well as the Hodge bound imply that
$$U_p(\Delta) \subseteq U_p(\Delta, T).$$
It follows that if $U_p(\Delta)$ is Zariski dense in $N(\Delta,
\mathbb{F}_p)$, then $U_p(\Delta, T)$ is also Zariski dense in
$N(\Delta, \mathbb{F}_p)$.

The Adolphson-Sperber conjecture \cite{AS} says that if $p\equiv
1~(\mod D)$, then $U_p(\Delta)$ is Zariski dense in $N(\Delta,
\mathbb{F}_p)$. This conjecture was proved to be true in
\cite{Wa1}\cite{Wa2} if $n\leq 3$. In particular, this implies

\begin{theorem} If $p\equiv 1~(\mod D)$ and $n\leq 3$, then
$U_p(\Delta, T)$ is Zariski dense in $N(\Delta, \mathbb{F}_p)$.
\end{theorem}

For $n\geq 4$, it was shown in \cite{Wa1}\cite{Wa2} that there is an
effectively computable positive integer $D^*(\Delta)$ depending only
on $\Delta$ such that if $p\equiv 1~(\mod D^*(\Delta))$, then
$U_p(\Delta)$ is Zariski dense in $N(\Delta, \mathbb{F}_p)$. Thus,
we obtain

\begin{theorem}\label{ord theorem} For each $\Delta$, there is an effectively computable positive
integer $D^*(\Delta)$ such that if $p\equiv 1~(\mod D^*(\Delta))$,
then  $U_p(\Delta, T)$ is Zariski dense in $N(\Delta,
\mathbb{F}_p)$.
\end{theorem}

The smallest possible $D^*(\Delta)$ is rather subtle to compute in
general, and it can be much larger than $D$. We now state a
conjecture giving reasonably precise information on $D^*(\Delta)$.

\begin{definition} Let
        $S(\triangle)$ be the monoid generated by the degree $1$
        lattice points in $M(\Delta)$, i.e., those lattice points on the
        codimension $1$ faces of $\Delta$ not containing the origin.
        Define the exponent of $\triangle$ by
        $$
        I(\triangle) = \inf\{ d\in \mathbb{Z}_{>0} | dM(\Delta) \subseteq S(\triangle)\}.$$
\end{definition}
If $u\in M(\Delta)$, then the degree of $Du$ will be integral but
$Du$ may not be a non-negative integral combination of degree $1$
elements in $M(\Delta)$ and thus $DM(\Delta)$ may not be a subset of
$S(\Delta)$. It is not hard to show that $I(\Delta)\geq D$. In
general they are different but they are equal if $n\leq 3$. This
explains why the Adolphson-Sperber conjecture is true if $n\leq 3$
and it can be false if $n\geq 4$. The following conjecture is a
modified form, and it is a consequence of Conjecture 9.1 in
\cite{Wa1}.

\begin{conjecture}\label{ord conjecture}
If $p \equiv 1~ mod ~ I(\triangle)$,  then  $U_p(\Delta)$ is Zariski
dense  in $N(\Delta, \mathbb{F}_p)$. In particular, $U_p(\Delta, T)$
is Zariski dense in $N(\Delta, \mathbb{F}_p)$ for such $p$.
\end{conjecture}

By
the facial decomposition theorem, in proving the above conjecture,
it is sufficient to assume that $\Delta$ has only one codimension
$1$ face not containing the origin.

\subsection{Generic Newton polygon} In the case that $U_p(\Delta,T)$
is empty, we expect the existence of a generic $T$-adic Newton
polygon. For this purpose, we need to re-scale the uniformizer. For
$f\in N(\Delta, \mathbb{F}_p)(\mathbb{F}_{p^a})$, the
$T^{a(p-1)}$-adic Newton polygon of $C_f(s,T;\mathbb{F}_{p^a})$ is
independent of the choice of $a$ for which $f$ is defined over
$\mathbb{F}_{p^a}$. We call them the absolute $T$-adic Newton
polygon of $f$.

\begin{conjecture}There is a Zariski open
dense subset $G_p(\Delta, T)$ of $N(\Delta, \mathbb{F}_p)$ such that
the absolute $T$-adic Newton polygon of $f$ is constant for all
$f\in G_p(\Delta,T)$. Denote this common polygon by ${\rm G
NP}_p(\Delta,T)$, and call it the generic Newton polygon of
$(\triangle,T)$.
\end{conjecture}

More generally, one expects that much of classical theory for finite
rank $F$-crystals extends to a certain nuclear infinite rank
setting. This includes the classical Dieudonne-Manin isogeny
theorem, the Grothendieck specialization theorem, the Katz isogeny
theorem \cite{Ka}. All these are essentially understood in the
ordinary infinite rank case, but open in the non-ordinary infinite
rank case.

Similarly, for each non-trivial $\psi$, there is a Zariski open
dense subset $G_p(\Delta, \psi)$ of $N(\Delta, \mathbb{F}_p)$ such
that the $\pi_{\psi}^{a(p-1)}$-adic Newton polygon of the $C$-value
$C_f(s,\pi_{\psi};\mathbb{F}_{p^a})$ is constant for all $f\in
G_p(\Delta,\psi)$. Denote this common polygon by ${\rm
GNP}_{p}(\Delta,\psi)$, and call it the generic Newton polygon of
$(\triangle,\psi)$. The existence of $G_p(\Delta, \psi)$ can be
proved. Since the non-degenerate assumption implies that the
C-function $C_f(s,\pi_{\psi})$ is determined by a single finite rank
$F$-crystal via a Dwork type cohomological formula for
$L_{f,\psi}(s)$. In the $T$-adic case, we are not aware of any such
finite rank reduction.

Clearly, we have the relation
$${\rm GNP}_{p}(\Delta,\psi) \geq {\rm GNP}_p(\Delta, T).$$
\begin{conjecture} If $p$ is sufficiently large, then
$${\rm GNP}_{p}(\Delta,\psi)={\rm GNP}_p(\Delta,T).$$
\end{conjecture}
This conjecture is proved in the case $n=1$ in \cite{LLN}.

Let ${\rm HP}(\Delta)$ denote the absolute Hodge polygon with
vertices $(0,0)$ and
$$(\sum\limits_{k=0}^iW(k), \sum\limits_{k=0}^i\frac{k}{D}W(k)),\ i=0,1,\cdots.
$$
Note that ${\rm HP}(\Delta)$ depends only on $\Delta$, not on $q$
any more. It is re-scaled from the $q$-Hodge polygon ${\rm
HP}_q(\Delta)$. Clearly, we have the relation
$${\rm GNP}_{p}(\Delta,\psi) \geq {\rm GNP}_p(\Delta, T) \geq {\rm HP}(\Delta).$$
 Conjecture \ref{ord
conjecture} says that if $p \equiv 1~(\mod ~ I(\triangle))$, then
both ${\rm GNP}_{p}(\Delta,\psi)$ and ${\rm GNP}_p(\Delta, T)$ are
equal to ${\rm HP}(\Delta)$. In general, the generic Newton polygon
lies above ${\rm HP}(\Delta)$ but for many $\Delta$ it should be
getting closer and closer to ${\rm HP}(\Delta)$ as $p$ goes to
infinity. We now make this more precise. Let $E(\Delta)$ be the
monoid generated by the lattice points in $\Delta$. This is a subset
of $M(\Delta)$. Generalizing the limiting Conjecture 1.11 in
\cite{Wa2} for $\psi$ of order $p$, we have
\begin{conjecture} If the difference $M(\Delta)-E(\Delta)$ is a
finite set, then for each non-trivial $\psi$, we have
$$\lim_{p\rightarrow \infty} {\rm GNP}_{p}(\Delta,\psi)={\rm HP}(\Delta).$$
In particular,
$$\lim_{p\rightarrow \infty}{\rm GNP}_p(\Delta,T)={\rm HP}(\Delta).$$
\end{conjecture}

This conjecture is equivalent to the existence of the limit. This is
because for all primes $p\equiv 1~(\mod D^*(\Delta))$, we already
have the equality ${\rm GNP}_{p}(\Delta,\psi)={\rm HP}(\Delta)$ by
Theorem \ref{ord theorem}. A stronger version of this conjecture
(namely, Conjecture 1.12 in \cite{Wa2}) has been proved by Zhu
\cite{Zh1}\cite{Zh2}\cite{Zhu3} in the case $m=1$ and $n=1$, see
also Blache and F\'erard \cite{BF}\cite{BFZ} and Liu \cite{Liu3} for
related further work in the case $m=1$ and $n=1$, Hong
\cite{H1}\cite{H2} and Yang \cite{Y} for more specialized one
variable results. For $n\geq 2$, the conjecture is clearly true for
any $\Delta$ for which both $D\leq 2$ and the Adolphson-Sperber
conjecture holds, because then ${\rm GNP}_p(\Delta,\psi)={\rm
HP}(\Delta)$ for every $p>2$. There are many such higher dimensional
examples \cite{Wa2}. Using free products of polytopes and the above
known examples, one can construct further examples \cite{B4}.

\subsection{$T$-adic Dwork Conjecture}

In this final subsection, we describe the $T$-adic version of Dwork's
conjecture \cite{Dw0} on pure slope zeta functions.

Let $\Lambda$ be a quasi-projective subvariety of $N(\Delta,
\mathbb{F}_p)$ defined over $\mathbb{F}_p$. Let $f_{\lambda}$ be a
family of Laurent polynomials parameterized by $\lambda \in
\Lambda$. For each closed point $\lambda \in \Lambda$, the Laurent
polynomial $f_{\lambda}$ is defined over the finite field
$\mathbb{F}_{p^{\deg(\lambda)}}$. The $T$-adic entire function
$C_{f_{\lambda}}(s,T)$ has the pure slope factorization
$$C_{f_{\lambda}}(s,T)=\prod_{\alpha \in \mathbb{Q}_{\geq 0}}
P_{\alpha}(f_{\lambda}, s),$$ where each $P_{\alpha}(f_{\lambda},
s)\in 1+s\mathbb{Z}_p[[T]][s]$ is a polynomial in $s$ whose
reciprocal roots all have $T^{{\deg(\lambda)}(p-1)}$-slope equal to $\alpha$.

\begin{definition}
For $\alpha \in \mathbb{Q}_{\geq 0}$, the $T$-adic pure slope
L-function of the family $f_{\Lambda}$ is defined to be the infinite
Euler product
$$L_{\alpha}(f_{\Lambda}, s) =\prod_{\lambda \in |\Lambda|}\frac{1}{
P_{\alpha}(f_{\lambda}, s^{\deg(\lambda)})}\in
1+s\mathbb{Z}_p[[T]][[s]],$$ where $|\Lambda|$ denotes the set of
closed points of $\Lambda$ over $\mathbb{F}_p$.
\end{definition}

The $T$-adic version of Dwork's conjecture is then the following

\begin{conjecture}
For $\alpha \in \mathbb{Q}_{\geq 0}$, the $T$-adic pure slope
L-function $L_{\alpha}(f_{\Lambda}, s)$ is $T$-adic meromorphic in
$s$.
\end{conjecture}

In the ordinary case, this conjecture can be proved using the
methods in \cite{Wa3}\cite{Wa4}\cite{Wa5}. It would be interesting
to prove this conjecture in the general case. The $\pi_{\psi}$-adic
version of this conjecture is essentially Dwork's original
conjecture, which can be proved as it reduces to finite rank
$F$-crystals. The difficulty of the $T$-adic version is that we have
to work with infinite rank objects, where much less is known in the
non-ordinary case.



\begin{thebibliography}{99}

\bibitem{AS}A. Adolphson and S. Sperber,
Exponential sums and Newton polyhedra: cohomology and estimates,
Ann. Math., 130 (1989), 367-406.
\bibitem{AS2}A. Adolphson and S. Sperber,
Newton polyhedra and the degree of the L-function associated to an
exponential sum, Invent. Math., 88(1987), 555-569.

\bibitem{BE}B. Berndt and R. Evans, The determination of Gauss sums, Bull. Amer. Math. Soc., 5(1981),
107-129.

\bibitem{Bl}R. Blache, Stickelberger's theorem for
$p$-adic Gauss sums, Acta Arith., 118(2005), no.1, 11-26.

\bibitem{BF}R. Blache and E. F\'erard, Newton straitification for
polynomials: the open stratum, J. Number Theory, 123(2007), 456-472.

\bibitem{BFZ}R. Blache, E. F\'erard and J.H. Zhu,
Hodge-Stickelberger polygons for L-functions of exponential sums of
$P(x^s)$, arXiv:0706.2340.

\bibitem{B4}R. Blache, Polygons de Newton de certaines sommes de
caract`eres et s'eries Poincar'e, arXiv:0802.3889.

\bibitem{B}E. Bombieri, On exponential sums in finite fields.
Amer. J. Math., 88(1966), 71-105.

\bibitem{Dw}B. Dwork, On the rationality of the zeta function of an
algebraic variety, Amer. J. Math., 82(1960), 631-648.

\bibitem{Dw0}B. Dwork, Normalized period matrices II, Ann. Math.,
98(1973), 1-57.

\bibitem{Dw1}B. Dwork, G. Gerotto and F.J. Sullivan, An Introduction
to G-Functions, Annals Math. Studies, Princeton University Press,
Number 133, 1994.



\bibitem{GKZ} I.M. Gelfand, M.M. Kapranov and A.V. Zelevinsky,
\it Discriminatns, Resultants and Multidimensional Determinants, \rm
Birkh\"user Boston, Inc., Boston, MA, 1994.

\bibitem{GK}B. Gross and N. Koblitz, Gauss sums and the $p$-adic $\Gamma$-functions,
Ann. Math., 109(1979), No. 2, 569-581.

\bibitem{Gr}A. Grothendieck, Formule de Lefschets et rationalit\'e
des fonctions L, S\'eminare Bourbaki, expos\'e 279, 1964/65.

\bibitem{H1}S. Hong, Newton polygons of L-functions associated with
exponential sums of polynomials of degree four over finite fields,
Finite Fields \& Appl., 7(2001), 205-237.

\bibitem{H2}S. Hong, Newton polygons of L-functions associated with
exponential sums of polynomials of degree six over finite fields, J.
Number Theory, 97(2002), 368-396.

\bibitem{Ig}J. Igusa, Forms of Higher Degree,
Tata Institute of Fundamental Research Lectures on Mathematics and Physics, 59,
by the Narosa Publishing House, New Delhi, 1978.

\bibitem{Ka}N. Katz, Slope filtration of F-crystals, Ast\'erisque,
63(1979), 113-164.

\bibitem{Liu}C.
Liu, The $L$-functions of twisted Witt coverings, J. Number Theory,
125 (2007), 267-284.
\bibitem{LW}C. Liu and D. Wei, The $L$-functions of Witt
coverings, Math. Z., 255 (2007), 95-115.

\bibitem{Liu3}C. Liu, Generic exponential sums associated to Laurent
polynomials in one variable, arXiv:0802.0271.

\bibitem{LLN}C. Liu, W. Liu, and C. Niu, Generic rigidity of Laurent
polynomials, arXiv:0901.0354.


\bibitem{Ma}Y. Manin, The theory of commutative formal groups over
fields of finite characteristic, Russian Math. Survey, 18(1963),
1-83.

\bibitem{Se}J-P. Serre, Endomorphismes compl\'etement continus des
espaces de Banach $p$-adiques, Publ. Math., IHES., 12(1962), 69-85.

\bibitem{Sp}S. Sperber, Congruence properties of hyperkloosterman
sums, Compositio Math., 40(1980), 3-33.

\bibitem{Wa1}D. Wan, Newton polygons of zeta functions and L-functions, Ann.
Math., 137 (1993), 247-293.
\bibitem{Wa2}D. Wan, Variation of $p$-adic Newton polygons for L-functions of exponential
sums, Asian J. Math., Vol 8, 3(2004), 427-474.

\bibitem{Wa3}D. Wan, Higher rank case of
Dwork's conjecture, J. Amer. Math. Soc., 13(2000), 807-852.

\bibitem{Wa4} D. Wan, Rank one case of Dwork's conjecture, J. Amer. Math. Soc.,
13(2000), 853-908.

\bibitem{Wa5}D. Wan, Dwork's conjecture on unit root zeta functions,
Ann. Math., 150(1999), 867-927.

\bibitem{Y}R. Yang, Newton polygons of L-functions of polynomials of
the form $x^d +\lambda x$, Finite Fields \& Appl., 9(2003), no.1,
59-88.


\bibitem{Zh1}J. H. Zhu, p-adic variation of L functions of one variable exponential sums,
I. Amer. J.  Math., 125 (2003), 669-690.

\bibitem{Zh2}J. H. Zhu, Asymptotic variation of L functions of one-variable exponential
sums, J. Reine Angew. Math., 572 (2004), 219--233.

\bibitem{Zhu3}J. H. Zhu, L-functions of exponential sums over one-dimensional affinoids :
Newton over Hodge, Inter. Math. Research Notices, no 30 (2004),
1529--1550.


\end{thebibliography}
\end{document}